\documentclass{amsart}

\usepackage{amssymb,latexsym}
\numberwithin{equation}{section}

\newtheorem{theorem}{Theorem}[section]
\newtheorem{proposition}[theorem]{Proposition}

\newtheorem{corollary}[theorem]{Corollary}
\newtheorem{remark}[theorem]{Remark}

\newtheorem{example}[theorem]{Example}
\newtheorem{definition}[theorem]{Definition}

\def\AA{\mathcal{A}}

\def\ZZ{\mathbb{Z}}
\def\CC{\mathbb{C}}
\def\RR{\mathbb{R}}

\def\BB{\mathcal{B}}
\def\MM{\mathcal{M}}
\def\TT{\mathbb{T}}
\def\gg{\mathfrak{g}}
\def\nn{\mathfrak{n}}
\def\hh{\mathfrak{h}}
\def\ii{\mathbf{i}}
\def\l{\ell}
\def\xx{\mathbf{x}}

\newcommand{\mat}[4]{\left(\!\!\begin{array}{cc}
#1 & #2 \\
#3 & #4 \\
\end{array}\!\!\right)}

\newcommand{\overunder}[2]{
\!\begin{array}{c}
\scriptstyle{#1}\\[-.1in]
-\!\!\!-\!\!\!-\\[-.1in]
\scriptstyle{#2}
\end{array}
\!
}

\begin{document}

\title[From LR-coefficients to cluster algebras]
{From Littlewood-Richardson coefficients to cluster algebras in
three lectures}

\author{Andrei Zelevinsky}
\address{\noindent Department of Mathematics, Northeastern University,
  Boston, MA 02115}
\email{andrei@neu.edu}

\date{\today}

\thanks{The author was supported in part
by NSF grant \#DMS-9971362.}

\maketitle

This is an expanded version of the notes of my three lectures at
a NATO Advanced Study Institute ``Symmetric functions 2001:
surveys of developments and perspectives"
(Isaac Newton Institute for Mathematical Sciences, Cambridge, UK;
June 25 -- July 6, 2001).
Lecture I presents a unified expression from \cite{bz01}
for generalized
Littlewood-Richardson coefficients (= tensor product
multiplicities) for any complex semisimple Lie algebra.
Lecture II outlines a proof of this result; the main idea of the
proof is to relate the LR-coefficients with canonical bases and
total positivity.
Lecture III introduces cluster algebras, a new class of
commutative algebras defined in \cite{fzclus1} in an attempt to
create an algebraic framework for canonical bases and total
positivity.

\tableofcontents

\section{Lecture I: Generalizing the Littlewood-Richardson rule}
\label{sec:lecture1}

\subsection{Schur polynomials and Littlewood-Richardson
coefficients}
Let $\Lambda_n = \ZZ[x_1, \dots, x_n]^{S_n}$ be the ring of
symmetric polynomials in $n$ variables.
Classical theory of symmetric polynomials studies
several linear bases in this ring.
All these bases are labeled by partitions
$\lambda = (\lambda_1 \geq \cdots \geq \lambda_n \geq 0)$ of length
$\leq n$.
The most important basis consists of \emph{Schur polynomials}
$s_\lambda = a_{\lambda + \delta}/a_\delta$, where
$a_\alpha = \det(x_i^{\alpha_j})$, and $\delta = (n-1, \dots, 1,0)$
(so that $a_\delta$ is the Vandermonde determinant
$\prod_{i < j} (x_i - x_j)$).
These polynomials were known long before Schur (they were
already known to Cauchy, and studied systematically by Jacobi) but I.~Schur was
the first to discover their representation-theoretic significance.
He established a bijective correspondence
$\lambda \mapsto V_\lambda$ between partitions of length $\leq n$
and polynomial representations of the group $GL_n(\CC)$.
Under this correspondence, the polynomial $s_\lambda$ is the character
of the corresponding representation $V_\lambda$:
$$s_\lambda (x_1, \dots, x_n) = {\rm tr} \ (g, V_\lambda) \ ,$$
where $x_1, \dots, x_n$ are the eigenvalues
(taken with multiplicities) of $g \in GL_n (\CC)$.

The \emph{Littlewood-Richardson coefficients} (LR-coefficients) are the structure
constants of the ring $\Lambda_n$ with respect to the basis of Schur
polynomials:
$$s_\lambda s_\nu = c_{\lambda \nu}^\mu s_\mu \ .$$
Since the product of characters corresponds to the tensor product
of representations, the LR-coefficients have the following
representation-theoretic interpretation: $c_{\lambda \nu}^\mu$
is the multiplicity of $V_\mu$ in the decomposition of
$V_\lambda \otimes V_\nu$ into the direct sum of irreducible
representations.
The \emph{Littlewood-Richardson rule} gives an explicit
combinatorial expression for $c_{\lambda \nu}^\mu$ in terms of
\emph{Young tableaux}.
This is one of the highlights of the classical theory.
Many different formulations and proofs of the rule are known by now
(see e.g., \cite{fomin-greene} and references there).
In this lecture, we present a more general result from \cite{bz01}
that puts the LR-rule into the general context of representations of
an arbitrary complex semisimple Lie algebra $\gg$.
Thus, we find an explicit combinatorial expression
for the multiplicity $c_{\lambda \nu}^\mu$ of $V_\mu$ in
$V_\lambda \otimes V_\nu$, where $V_\lambda$ is an irreducible
finite-dimensional representation of $\gg$ with highest weight $\lambda$.
(This expression is different from the one
obtained by P.~Littelmann in \cite{lit} in terms of his path model.)

\subsection{Semisimple Lie algebras}
Let us briefly review some basic definitions and
properties of complex semisimple Lie algebras.

\begin{definition}
\label{def:Cartan-matrix}
{\rm A \emph{Cartan matrix} of rank $r$ is an $r \times r$ integer
matrix $A = (a_{ij})$ such that

\noindent (1) $a_{ii} = 2$ for all $i$;

\noindent (2) $a_{ij} \leq 0$ for $i \neq j$, and
$a_{ij} = 0  \Rightarrow a_{ji} = 0$;

\noindent (3) for every non-empty subset
$I \subset [1,r] = \{1,2, \dots, r\}$, the principal
minor $\det ((a_{ij})_{i,j \in I})$ is positive.}
\end{definition}

The list of all Cartan matrices is given by the famous Cartan-Killing classification.
It asserts that every Cartan matrix can be transformed by a simultaneous
permutation of rows and columns into a direct sum of matrices of
types $A_r (r \geq 1), B_r (r \geq 2), C_r (r \geq 3),
D_r (r \geq 4), E_6, E_7, E_8, F_4, G_2$.
Here the subscript indicates the rank.
Our standard example will be the type $A_r$ for which $a_{ij} = -1$
if $|i - j| = 1$, and $a_{ij} = 0$ if $|i - j| > 1$.

\begin{definition}
\label{def:g(A)}
{\rm The complex semisimple Lie algebra $\gg = \gg (A)$ associated to a
Cartan matrix $A$ is generated by $3r$ Chevalley generators
$e_i, \alpha_i^\vee$, and $f_i$ for $i=1, \ldots, r$ with the defining
relations
$$[\alpha_i^\vee, \alpha_j^\vee] = 0, \,\,  [\alpha_i^\vee, e_j] =
a_{ij} e_j, \,\,  [\alpha_i^\vee, f_j] = - a_{ij} f_j,
\,\, [e_i, f_j] = \delta_{ij} \alpha_i^\vee \ ;$$
$$({\rm ad} \ e_i)^{1 - a_{ij}} (e_j) = ({\rm ad} \ f_i)^{1 - a_{ij}} (f_j)
= 0 \,\, (i \neq j) \ .$$
}
\end{definition}

Note that the definition of $\gg (A)$ makes sense for every
integer matrix $A$ satisfying conditions (1) and (2) in
Definition~\ref{def:Cartan-matrix} (this is essentially a
Kac-Moody algebra); condition (3) is necessary and sufficient for
$\gg(A)$ to be finite-dimensional.

The elements $\alpha_i^\vee$ are the \emph{simple coroots} of $\gg$;
they form a basis of an abelian subalgebra $\hh$ of $\gg$ called
the Cartan subalgebra. The elements $e_i$ (resp.,~$f_i$) are called
\emph{raising} (resp.,~\emph{lowering}) generators; they generate
a maximal nilpotent Lie subalgebra $\nn$ (resp.,~$\nn_-$) of
$\gg$. The algebra $\gg$ has the \emph{triangular decomposition}
$\gg = \nn_- \oplus \hh \oplus \nn$. The \emph{simple roots}
$\alpha_1, \dots, \alpha_r$ form a basis in the dual space $\hh^*$
such that $[h, e_i] = \alpha_i (h) e_i$, and $[h,f_i] = - \alpha_i
(h) f_i$ for any $h \in \hh$ and $i \in [1,r]$. Thus, the Cartan
matrix $A = (a_{ij})$ is given by
$a_{ij} = \langle \alpha_i^\vee, \alpha_j\rangle$, where
$\langle h, \alpha\rangle$ stands for the natural pairing
$\hh \times \hh^* \to \CC$.

The \emph{Weyl group} $W$ of $\gg$ is a group of linear transformations
of $\hh^*$ generated by \emph{simple reflections}
$s_1, \dots, s_r$ acting by
$s_i (\gamma) = \gamma - \langle \alpha_i^\vee, \gamma\rangle \alpha_i$.
As an abstract group, $W$ is a finite Coxeter group; we denote by
$\l (w)$ the length function on $W$, i.e., the minimal length of
a factorization of $w$ into a product of simple reflections.
A \emph{reduced word} for $w \in W$ is a sequence of indices
$(i_1, \ldots, i_l)$ that satisfies $w = s_{i_1} \cdots s_{i_l}$
and has the shortest possible length~$l=\l(w)$. The set of all reduced
words for~$w$ will be denoted by~$R(w)$. As is customary, $w_0$
denotes the unique element of maximal length in~$W$.

The \emph{weight lattice} $P$ of $\gg$ consists of all $\gamma \in \hh^*$ such that
$\langle \alpha_i^\vee, \gamma \rangle \in \ZZ$ for all $i$.
Thus $P$ has a $\ZZ$-basis $\omega_1, \dots, \omega_r$ of \emph{fundamental weights}
given by $\langle \alpha_i^\vee, \omega_j \rangle = \delta_{i,j}$.
A weight $\lambda \in P$ is \emph{dominant} if
$\langle \alpha_i^\vee, \lambda\rangle \geq 0$ for
any $i \in [1,r]$; thus $\lambda = \sum_i l_i \omega_i$ with all $l_i$ nonnegative integers.
The standard partial order on $P$ is given as follows: $\beta < \alpha$
means that $\alpha - \beta $ is a positive integer linear
combination of simple roots.

\begin{example}
\label{ex:A_r}
{\rm For the type $A_r$, the Lie algebra $\gg$ is $sl_{r+1}$,
the algebra of $(r+1) \times (r+1)$-matrices with trace $0$.
The Chevalley generators are
$e_i = E_{i,i+1}$, $\alpha_i^\vee = E_{i,i} - E_{i+1,i+1}$,
$f_i = E_{i+1,i}$, where the $E_{i,j}$ are matrix units.
The Cartan subalgebra $\hh$ consists of diagonal matrices
with trace~$0$, while $\nn$ (resp.~$\nn_-$) consists of nilpotent
upper-triangular (resp.~lower-triangular) matrices.
Let $\varepsilon_1, \dots, \varepsilon_{r+1}$ denote the linear
forms on the space of diagonal matrices given by
$\langle E_{i,i}, \varepsilon_j \rangle = \delta_{i,j}$.
The weight lattice $P$ is canonically identified with
$(\ZZ \varepsilon_1 \oplus \cdots \oplus \ZZ \varepsilon_{r+1})/
\ZZ (\varepsilon_1 + \cdots + \varepsilon_{r+1})$.
The simple roots are $\alpha_i = \varepsilon_i - \varepsilon_{i+1}$.
The Weyl group $W$ is the symmetric group $S_{r+1}$ acting on $P$
by permuting the forms $\varepsilon_i$.
The simple reflection $s_i$ is the transposition $(i,i+1)$.
The fundamental weights are
$\omega_i = \varepsilon_1 + \cdots + \varepsilon_{i}$.
}
\end{example}

Let $\gg^\vee$ denote the complex semisimple Lie algebra
associated as above with the transpose Cartan matrix $A^T$
(this algebra is sometimes referred to as the \emph{Langlands dual} of $\gg$).
Thus we can identify the Cartan subalgebra of $\gg^\vee$ with $\hh^*$, and the simple roots
(resp., coroots) of $\gg^\vee$ with simple coroots (resp., roots) of $\gg$.
We denote the fundamental weights for $\gg^\vee$ by $\omega_1^\vee, \dots, \omega_r^\vee$;
they are the elements of $\hh$ such that $\langle \omega_i^\vee, \gamma\rangle$ is the coefficient
of $\alpha_i$ in the expansion of $\gamma \in \hh^*$ in the basis of simple roots.
The Weyl groups of $\gg$ and $\gg^\vee$ are naturally identified with each other, and we denote
both groups by the same symbol $W$.
Thus, the simple reflections act in $\hh$ by
$s_i (h) = h - \langle h, \alpha_i\rangle \alpha_i^\vee$.

Every finite-dimensional $\gg$-module $V$ is known to be $\hh$-diagonalizable,
i.e., it has the \emph{weight decomposition} $V = \oplus_{\gamma \in P}  V(\gamma)$,
where $V(\gamma) = \{v \in V: h v = \gamma (h) v \, \text {\, for \,} h \in \hh\}$.
Note that $e_i (V(\gamma)) \subset V(\gamma + \alpha_i)$ and
$f_i (V(\gamma)) \subset V(\gamma - \alpha_i)$, which is why we
call the $e_i$ (resp.,~$f_i$) the raising (resp.,~lowering) generators.
Every simple finite-dimensional $\gg$-module has a unique \emph{highest
weight} $\lambda$ (i.e., the maximal weight with respect to the
standard partial order), and this correspondence provides a
bijection $\lambda \mapsto V_\lambda$ between dominant weights
and isomorphism classes of simple finite-dimensional $\gg$-modules.
Extending the LR-coefficients,
we denote by $c_{\lambda, \nu}^\mu$ the multiplicity of $V_\mu$
in the tensor product $V_\lambda \otimes V_\nu$.
Our aim is to give an explicit expression for these multiplicities.

\subsection{Semisimple Lie groups}
Let $G$ be a simply connected connected complex semisimple Lie group with the
Lie algebra $\gg$.
For $i \in [1,r]$, we denote by $x_i (t)$ and $y_i (t), \, t \in \CC$, the one-parameter subgroups
in $G$ given by
$$x_i (t) = \exp \ (t e_i), \,\, y_i (t) = \exp \ (t f_i) \ .$$
Let $N$ (resp., $N_-$) be the maximal unipotent subgroup of $G$ generated by all
$x_i (t)$ (resp., $y_i (t)$).
Let $H$ be the maximal torus in $G$ with the Lie algebra~$\hh$.
Let $B = HN$ and $B_- = HN_-$ be two
opposite Borel subgroups, so that $H=B_-\cap B$.
For every $i\in [1,r]$, let $\varphi_i: SL_2 \to G$ denote
the canonical embedding corresponding to the simple
root~$\alpha_i\,$; thus we have
$$x_i (t) = \varphi_i \mat{1}{t}{0}{1}, \,\, y_i (t) = \varphi_i \mat{1}{0}{t}{1} \ .$$
We also set
$$t^{\alpha_i^\vee} = \varphi_i \mat{t}{0}{0}{t^{-1}} \in H$$
for any $i$ and any $t \neq 0$.

The Weyl group $W$ of $\gg$ is naturally identified with ${\rm Norm}_G (H)/H$;
this identification sends each simple reflection
$s_i$ to the coset $\overline {s_i} H$, where the representative
$\overline {s_i} \in {\rm Norm}_G (H)$ is defined by
\[
\overline {s_i} = \varphi_i \mat{0}{-1}{1}{0} \, .
\]
The elements $\overline {s_i}$ satisfy the braid relations in~$W$;
thus the representative $\overline w$ can be unambiguously defined for any
$w \in W$ by requiring that
$\overline {uv} = \overline {u} \cdot \overline {v}$
whenever $\l (uv) = \l (u) + \l (v)$.

The weight lattice $P$ is identified with the group of
multiplicative characters of~$H$,
written here in the exponential notation:
a weight $\gamma\in P$ acts by $h \mapsto h^\gamma$.
Under this identification, the fundamental weights
$\omega_1, \ldots, \omega_r$ are given by
$(t^{\alpha_i^\vee})^{\omega_j} = t^{\delta_{i,j}}$.
The action of~$W$ on $P$ can be now written as
$h^{w (\gamma)} = (w^{-1} h w)^\gamma$
for $w \in W$, $h \in H$, $\gamma \in P$.

\subsection{A geometric realization of LR-coefficients}
For an algebraic variety $X$, let $\CC[X]$ denote the algebra of
regular (algebraic) functions on $X$.
In particular, we consider the algebra $\CC[N]$.
Since $H$ normalizes $N$, the algebra $\CC[N]$
has a $Q_+$-grading, where $Q_+$
is the additive subsemigroup in $\hh^*$ spanned by the simple
roots $\alpha_i$: for $\gamma \in Q_+$, the corresponding
homogeneous component is given by
$$\CC[N](\gamma) = \{F \in \CC[N]: F(h n h^{-1}) = h^\gamma
F(n) \text {\, for \,} h \in H, n \in N\} \ .$$
Let us define linear operators $L_i$ and $R_i$ on
$\CC[N]$ by
$$L_i F (n) = \frac{d}{dt} \Biggl |_{t=0} F(x_i (t) n), \, \,
R_i F (n) = \frac{d}{dt} \Biggl |_{t=0} F(n x_i (t)) \ .$$

\begin{proposition}
\label{pr:c-on-N}
For any three dominant weights $\lambda$, $\mu$ and $\nu$, we have
\begin{equation*}
c_{\lambda, \nu}^\mu
= \dim \ \{F \in \CC[N](\lambda + \nu - \mu):
L_i^{\langle \alpha_i^\vee, \lambda\rangle + 1}F =
R_i^{\langle \alpha_i^\vee, \nu\rangle + 1}F = 0
\text {\,\, for \,\,} i \in [1,r]\} \ .
\end{equation*}
\end{proposition}

This proposition is classical.
Let us outline the proof, which proceeds in three steps.
The first step is to use the following
interpretation of multiplicities (see~\cite{prv}):
\begin{equation}
\label{eq:prv}
c_{\lambda, \nu}^\mu = \dim V_\lambda (\mu - \nu; \nu) \ ,
\end{equation}
where
\begin{equation}
\label{eq:prv-space}
V_\lambda (\gamma; \nu) =
\{v \in V_\lambda (\gamma):
e_i^{\langle \alpha_i^\vee, \nu\rangle + 1} (v) = 0
\text  {\,\, for  all \,\,}  i\} \ .
\end{equation}

Second, for every dominant weight $\lambda$, we shall denote by the same
symbol $V_\lambda$ the irreducible finite-dimensional
representation of $G$ obtained from the $\gg$-module $V_\lambda$
via exponentiation.
An explicit geometric realization of $V_\lambda$ can be given as
follows.
Note that if $G$ acts on the right in an algebraic variety
$X$, then $\CC[X]$ becomes a linear
representation of $G$ via $gF (x) = F(xg)$.
We apply this construction to $X = N_- \backslash G$; this variety is called
the \emph{base affine space} of $G$.
We have
$$\CC[N_- \backslash G] = \{F \in \CC[G]: F(n_- g) = F(g) \text {\,\, for \,\,}
g \in  G, n_-  \in N_-\} \ .$$
Since $H$ normalizes $N_-$, it acts on the left
in $N_- \backslash G$ via $h(N_- g) = hN_-g = N_-hg$.
Then $V_\lambda$ can be realized as a
subrepresentation in $\CC[N_- \backslash G]$ given by
$$V_\lambda = \{F \in \CC[N_- \backslash G]: F(h x) = h^\lambda F(x)
\text {\, for \,} x \in N_- \backslash G, h \in H\} \ .$$
Furthermore, $\CC[N_- \backslash G]$ is the direct sum of all $V_\lambda$,
for all dominant weights~$\lambda$.

Finally, the above realization can be modified as follows.
Let $G_0=N_-HN$ denote the open subset of elements $x\in G$ that
have a Gaussian decomposition; this (unique) decomposition will be written as
$x = [x]_- [x]_0 [x]_+ \,$.
Thus, we have an open embedding
$B = N_- \backslash G_0 \hookrightarrow N_- \backslash G$
inducing an embedding (via restriction)
$\CC[N_- \backslash G] \hookrightarrow \CC[B]$.
Each $F \in V_\lambda$ viewed as a function on $B$ satisfies
$F(hn) = h^\lambda F(n)$ for $n \in N$ and $h \in H$.
Thus, $F$ is uniquely determined by its restriction to $N$,
which allows us to embed $V_\lambda$ into $\CC[N]$.
(Note that in this realization, different irreducible
representations $V_\lambda$ are no longer disjoint.)
It remains to notice that, under this embedding,
the subspace $V_\lambda (\gamma; \nu)$
(see (\ref{eq:prv-space})) becomes the following subspace in $\CC[N]$:
\begin{equation}
\label{eq:prv-concrete}
V_\lambda (\gamma; \nu)
= \{F \in \CC[N](\lambda  - \gamma):
L_i^{\langle \alpha_i^\vee, \lambda\rangle + 1}F =
R_i^{\langle \alpha_i^\vee, \nu\rangle + 1}F = 0
\text  {\,\, for  all \,\,}  i\} \ .
\end{equation}

\subsection{Generalized minors}
Following \cite{fz}, we now introduce a family of regular
functions on $G$ generalizing the minors of a matrix.
For $u,v \in W$ and $i \in [1,r]$, the \emph{generalized minor}
$\Delta_{u \omega_i, v \omega_i}$
is the regular function on $G$ whose restriction to the open set
${\overline {u}} G_0 {\overline {v}}^{-1}$ is given by
\begin{equation}
\label{eq:Delta-general}
\Delta_{u \omega_i, v \omega_i} (x) =
(\left[{\overline {u}}^{\ -1}
   x \overline v\right]_0)^{\omega_i} \ .
\end{equation}
As shown in \cite{fz}, $\Delta_{u \omega_i, v \omega_i}$ is well defined, and
it depends on the weights $u \omega_i$ and $v \omega_i$ alone, not on the particular
choice of $u$ and~$v$.
In the special case of type $A_r$, where $G=SL_{r+1}$,
we have $u \omega_i =  \sum_{i \in I}\varepsilon_i$ and
$v \omega_i =  \sum_{i \in J}\varepsilon_i$ for some subsets
$I, J \subset [1,r+1]$ of the same size $i$; it is then an easy exercise
to show that $\Delta_{u \omega_i, v \omega_i} = \Delta_{I,J}$, the
minor with the row set $I$ and the column set $J$.

Many familiar properties of minors can be extended to generalized minors.
For example, we have
$\Delta_{\gamma, \delta} (x^T) = \Delta_{\delta, \gamma} (x)$,
where
$x \mapsto x^T$ is the ``transpose" involutive antiautomorphism of
$G$ given by
\begin{equation}
\label{eq:T}
h^T = h \quad (h \in H) \ , \quad x_i (t)^T = y_i (t) \ ,
\quad y_i (t)^T = x_i (t) \ .
\end{equation}

Now we present some less obvious identities for generalized minors.
The following identity was obtained in \cite[Theorem~1.17]{fz}.
For the type $A_r$, it plays a crucial role in C.~L.~Dodgson's condensation method,
and is because of that occasionally associated with the name of Lewis
Carroll.

\begin{proposition}
\label{pr:minors-Dodgson}
Suppose $u,v \in W$ and $i \in [1,r]$
are such that $\l (us_i) = \l (u) + 1$ and $\l (vs_i) = \l (v) + 1$.
Then
\[
\Delta_{u \omega_i, v \omega_i} \Delta_{us_i \omega_i, v s_i \omega_i}
= \Delta_{us_i \omega_i, v \omega_i} \Delta_{u \omega_i, v s_i \omega_i}
+ \prod_{j \neq i} \Delta_{u \omega_j, v \omega_j}^{- a_{ji}} \ .
\]
\end{proposition}

The next proposition presents a special case of generalized Pl\"ucker relations
first obtained in~\cite[Corollary~6.6]{bz97}.

\begin{proposition}
\label{pr:minors-Plucker}
Let $w \in W$ and $i, j \in [1,r]$ be such that
$a_{ij}= a_{ji} = -1$ and $\l (ws_i s_j s_i)=\l(w) + 3$.
Then
\[
\Delta_{\omega_i,ws_i \omega_i} \Delta_{\omega_j, ws_j\omega_j} =
\Delta_{\omega_i, w \omega_i} \Delta_{\omega_j, ws_i s_j\omega_j} +
\Delta_{\omega_i, w s_j s_i \omega_i} \Delta_{\omega_j, w\omega_j} \ .
\]
\end{proposition}

If $a_{ij} a_{ji}$ is equal to $2$ or $3$, the
corresponding Pl\"ucker relations are more complicated;
see \cite{bz97,bz01,fz}.

We will be especially interested in the restrictions of
generalized minors to $B$ or $N$.
One can show that $\Delta_{\gamma, \delta}\Bigl |_{B} = 0$ unless
$\gamma \geq \delta$.
More precisely, the following characterization follows easily
from \cite[Theorem~5.8]{bz01}.

\begin{proposition}
\label{pr:B-minors}
Let $\gamma = u \omega_i$ and $\delta = v \omega_i$.
The minor $\Delta_{\gamma, \delta}$ has non-zero
restriction to $B$ if and only if the weight subspace
$V_{\omega_i} (\gamma)$ is contained in the
$\nn$-submodule of $V_{\omega_i}$ generated by
$V_{\omega_i} (\delta)$.
\end{proposition}

We shall call the pairs $(\gamma, \delta)$ as in
Proposition~\ref{pr:B-minors} (as well as the corresponding minors)
\emph{upper-triangular}.

\subsection{Tropicalization and the main result}
\label{sec:main}
Let us introduce a family of commuting variables
$M_{\gamma, \delta}$ labeled by
upper-triangular pairs of weights $(\gamma, \delta)$.
Let $\MM$ denote the variety of tuples
$(M_{\gamma, \delta})$, whose components satisfy
the relations in Propositions~\ref{pr:minors-Dodgson} and
\ref{pr:minors-Plucker}, where all non-upper-triangular minors are
specialized to $0$, and each $\Delta_{\gamma, \delta}$ is
replaced by $M_{\gamma,\delta}$
(in the non-simply-laced case, we also impose the remaining Pl\"ucker relations
mentioned above).
The correspondence
$b \mapsto (M_{\gamma,\delta} = \Delta_{\gamma,\delta}(b))$
is a well-defined regular map $B \to \MM$;
the results in \cite{fz} imply that it is a birational isomorphism.

An important feature of $\MM$ is that all the defining relations
are \emph{subtraction-free}, i.e., they only involve addition and multiplication.
This makes it possible to consider the \emph{tropical specialization}
$\MM_{\rm trop}$, taking all components $M_{\gamma, \delta}$ to be integers, and
replacing the usual addition and multiplication
with their tropical versions given by
\begin{equation}
\label{eq:min-+}
a \,\oplus \,b = {\rm min}\, (a,b) \ , \quad a\, \odot \, b = a + b \ .
\end{equation}
We shall also denote by $\MM^\vee$ and $\MM_{\rm trop}^\vee$ the
corresponding varieties associated with $\gg^\vee$, the Langlands dual Lie algebra.

Finally, everything is ready for presenting our expression for the
generalized LR-coefficients (cf.~\cite[Theorem~5.15]{bz01}).

\begin{theorem}
\label{th:LR-Plucker-Lusztig}
For any three dominant weights $\lambda, \mu, \nu$ for $\gg$, the multiplicity
$c_{\lambda, \nu}^\mu$ is equal to the number of integer tuples
$(M_{\gamma,\delta}) \in  \MM_{\rm trop}^\vee$
satisfying the following conditions for any $i \in [1,r]$:

\smallskip

\noindent {\rm (0)} $M_{\omega_i^\vee, \omega_i^\vee} = 0$;

\smallskip

\noindent {\rm (1)} $M_{\omega_i^\vee, s_i \omega_i^\vee} \geq 0$;

\smallskip

\noindent {\rm (2)} $M_{\omega_i^\vee, w_0 \omega_i^\vee} =
\langle \omega_i^\vee, \lambda + \nu - \mu \rangle$;

\smallskip

\noindent {\rm (3)} $M_{s_i \omega_i^\vee, w_0 \omega_i^\vee}  \geq
\langle \omega_i^\vee, s_i \lambda + \nu - \mu \rangle$;

\smallskip

\noindent {\rm (4)} $M_{\omega_i^\vee, w_0 s_i\omega_i^\vee} \geq
\langle \omega_i^\vee, \lambda +  s_i \nu - \mu \rangle$.
\end{theorem}

It would be interesting to find a geometric proof of
Theorem~\ref{th:LR-Plucker-Lusztig} in the spirit of \cite{bg,mv}.
An outline of the proof from \cite{bz01} will be presented in the next lecture.

\section{Lecture II: Canonical bases and total positivity}
\label{sec:lecture2}

\subsection{Dual realization of $c_{\lambda, \nu}^\mu$}
\label{sec:dual}
The universal enveloping algebra $U(\nn)$
is $Q_+$-graded via ${\rm deg} (e_i) = \alpha_i$ for all $i$.
Realizing elements of $U(\nn)$ as left-invariant differential operators on $N$,
we obtain, for every $\gamma \in Q_+$, a natural non-degenerate pairing between
$U(\nn)(\gamma)$ and $\CC[N](\gamma)$, thus identifying $U(\nn)(\gamma)$
with the dual space $(\CC[N](\gamma))^*$.
Under this identification, the dual of the subspace
$V_\lambda (\gamma; \nu)$ defined in
(\ref{eq:prv-concrete}) becomes the homogeneous component of degree
$\lambda - \gamma$ in
$$U(\nn) \Bigl / \sum_i \left (e_i^{\langle \alpha_i^\vee, \lambda \rangle + 1} U(\nn) \ + \
U(\nn) e_i^{\langle \alpha_i^\vee, \nu \rangle + 1}\right ) \ .$$
Therefore, to prove Theorem~\ref{th:LR-Plucker-Lusztig}
it suffices to show the following.

\begin{proposition}
\label{pr:Un-basis}
There exists a linear basis $\BB$ in $U(\nn)$ satisfying the
following properties:

\begin{itemize}

\item[{\rm (1)}]
$\BB$ consists of homogeneous elements, and
every subspace of the form
$$\sum_i (e_i^{l_i + 1} U(\nn) \ + \  U(\nn) e_i^{n_i+ 1})
\,\, (l_i, n_i \in \ZZ_{\geq 0})$$
is spanned by a part of $\BB$.

\item[{\rm (2)}] $\BB$ is labeled by the set of
integer tuples $(M_{\gamma,\delta}) \in  \MM_{\rm trop}^\vee$ such that
$M_{\omega_i^\vee,\omega_i^\vee} = 0$ and
$M_{\omega_i^\vee, s_i \omega_i^\vee} \geq 0$
for all $i$.

\item[{\rm (3)}] An element of $\BB$ has degree $\lambda + \nu - \mu$ and
does not belong to the subspace
$\sum_i (e_i^{\langle \alpha_i^\vee, \lambda \rangle + 1} U(\nn) \ + \
U(\nn) e_i^{\langle \alpha_i^\vee, \nu \rangle + 1})$
if and only if the corresponding tuple
$(M_{\gamma,\delta})$ satisfies conditions {\rm (2)--(4)} in
Theorem~\ref{th:LR-Plucker-Lusztig}.

\end{itemize}
\end{proposition}

We shall show that $\BB$ can be chosen as
the specialization at $q = 1$ of Lusztig's \emph{canonical basis}
in the quantized universal enveloping algebra $U_q (\nn)$.

\subsection{Canonical bases and their Lusztig parametrizations}
\label{sec:canonical bases}
Let us recall some basic facts about quantized universal enveloping algebras and
their canonical bases.
Unless otherwise stated, all results in this section are due to G.~Lusztig and
can be found in \cite{lu}.
The quantized universal enveloping
algebra $U=U_q (\gg)$ associated to $\gg$ is defined as follows.
Fix positive integers $d_1, \dots, d_r$ such that $d_i a_{ij}=d_j a_{ji}$,
where $(a_{ij})$ is the Cartan matrix of $\gg$.
The algebra $U$ is a $\CC(q)$-algebra with unit generated by the elements
$E_i, K_i^{\pm 1}$, and $F_i$ for $i = 1, \dots, r$
subject to the relations
$$K_i K_j = K_j K_i, \,\, K_i E_j K_i^{-1} = q^{d_i a_{ij}}E_j, \,\,
K_iF_jK_i^{-1}=q^{-d_ia_{ij}}E_j \ ,$$
$$E_i F_j - F_j E_i= \delta_{ij} \frac{K_i-K_i^{-1}}{q^{d_i}-q^{-d_i}}$$
for all $i$ and $j$, and the \emph{quantum Serre relations}
$$\sum_{k+l=1-a_{ij}} (-1)^k  E_i^{(k)} E_j E_i^{(l)} = \sum_{k+l=1-a_{ij}}
(-1)^k F_i^{(k)}F_jF_i^{(l)} = 0$$
for $i\ne j$.
Here $E_i^{(k)}$ and $F_i^{(k)}$ stand for the \emph{divided powers}
defined by
$$E_i^{(k)}=\frac{1}{[1]_i [2]_i \cdots [k]_i}E_i^k, \,\,
F_i^{(k)}=\frac{1}{[1]_i [2]_i \cdots [k]_i} F_i^k \ ,$$
where $[l]_i=\frac{q^{d_i l}-q^{-d_i l}}{q^{d_i }-q^{-d_i }}$.
The algebra $U$ is graded by the root lattice of $\gg$ via
$${\rm deg}(K_i)=0, \,\, {\rm deg}(E_i)= -{\rm deg}(F_i)= \alpha_i \ .$$

To each $i = 1, \dots, r$, Lusztig associates an
algebra automorphism $T_i$ of $U$ uniquely determined by:
$$T_i (K_j)= K_j K_i^{- a_{ij}} \quad (j = 1, \dots, r) \ ,$$
$$T_i (E_i)= -K_i^{-1} F_i, \,\, T_i(F_i)= -E_i K_i \ ,$$
and, for all $j\neq i$,
$$T_i (E_j)= \sum_{k+l= -a_{ij}} (-1)^k q^{- d_i k} E_i^{(k)}E_j E_i^{(l)}
\ ,$$
$$T_i (F_j)=\sum_{k+l=-a_{ij}} (-1)^k q^{d_i k} F_i^{(l)}F_jF_i^{(k)} \,\,
\ .$$
(This automorphism was denoted by $T'_{i,-1}$ in \cite{lu}.)
The $T_i$ satisfy the braid relations
and so extend to an action of the braid group on $U$.

Let $U^+ = U_q(\nn)$ denote the subalgebra of $U$ generated by $E_1, \dots, E_r$.
We now recall Lusztig's definitions
of the PBW-type bases and the canonical basis in $U^+$.
For a reduced word $\ii=(i_1,\ldots,i_m) \in R(w_0)$, and
an $m$-tuple $t = (t_1, \dots, t_m) \in \ZZ^m_{\ge 0}$, denote
$$p_\ii^{(t)} := E_{i_1}^{(t_1)} T_{i_1}(E_{i_2}^{(t_2)})
\cdots (T_{i_1} \cdots T_{i_{m-1}})(E_{i_m}^{(t_m)}) \ .$$
As shown in \cite{lu}, all these elements belong to $U^+$.
For a given $\ii$, the set of all $p_\ii^{(t)}$ with $t \in \ZZ^m_{\ge 0}$ is called
the \emph{PBW type basis} corresponding to $\ii$ and is denoted by
$\BB_\ii$.
This terminology is justified by the following proposition
proved in \cite[Corollary~40.2.2]{lu}.

\begin{proposition}
For every $\ii \in R(w_0)$, the set $\BB_\ii$ is a $\CC(q)$-basis of $U_+$.
\end{proposition}

According to \cite[Proposition~8.2]{lu96}, the canonical basis $\BB$ of $U^+$ can
now be defined as follows.
Let $u \mapsto \overline u$ denote the
$\CC$-linear involutive algebra automorphism of $U^+$ such that
$\overline {q} = q^{-1}, \,\, \overline {E_i} = E_i$.

\begin{proposition}
\label{pr:characterization of BB}
For every $\ii \in R(w_0)$ and $t \in  \ZZ_{\ge 0}^m$,
there is a unique element $b = b_\ii (t)$ of $U^+$ such that
$\overline {b} = b$, and $b - p_\ii^{(t)}$ is a linear combination
of the elements of $\BB_\ii$ with coefficients in $q^{-1} \ZZ[q^{-1}]$.
The set $\{b_{\ii}(t): t \in  \ZZ_{\ge 0}^m\}$ does not depend on the choice
of $\ii$; it is called the canonical basis and denoted by $\BB$.
\end{proposition}

In view of Proposition~\ref{pr:characterization of BB},
any $\ii \in R(w_0)$ gives rise to a bijection
$b_\ii: \ZZ_{\ge 0}^m \to \BB$.
We refer to these bijections as \emph{Lusztig parametrizations} of $\BB$.
Let us summarize some of their properties.
To do this, we need some more notation.

Let $i \mapsto i^*$ denote the involution on $[1,r]$
defined by $w_0 (\alpha_i) = - \alpha_{i^*}$.
For every sequence $\ii = (i_1, \dots, i_m)$,
we denote by $\ii^*$ and $\ii^{\rm op}$ the sequences
\begin{equation}
\label{eq:* and op}
\ii^* = (i_1^*, \dots, i_m^*), \,\, \ii^{\rm op} = (i_m, \dots, i_1) \ ;
\end{equation}
clearly, both operations $\ii \mapsto \ii^*$ and $\ii \mapsto \ii^{\rm op}$
preserve the set of reduced words $R(w_0)$.

\begin{proposition}
\label{pr:Lusztig parametrizations}

\noindent
\begin{itemize}
\item[{\rm (i)}] Any canonical basis vector $b_\ii (t) \in \BB$ is homogeneous of degree
$\sum_k t_k \cdot s_{i_1} \cdots s_{i_{k-1}} \alpha_{i_k}$.

\item[{\rm (ii)}] Every subspace of the form $E_i^n U^+$ in $U^+$ is spanned by
a subset of $\BB$.
Furthermore, let $l_i (b)$ denote the maximal integer $n$ such that $b \in E_i^n U^+$;
then, for any $\ii \in R(w_0)$ which begins with $i_1 = i$, we have
$l_i (b_\ii (t_1, \dots, t_m)) = t_1$.

\item [{\rm (iii)}] The canonical basis $\BB$ is stable under
the involutive $\CC(q)$-linear algebra antiautomorphism $E \to E^\iota$ of
$U^+$ such that $E_i^\iota = E_i$ for all $i$.
Furthermore, we have $b_\ii (t)^\iota = b_{\ii^{* {\rm op}}} (t^{\rm op})$.
\end{itemize}
\end{proposition}

As a corollary, we obtain the following interpretation of
$c_{\lambda, \nu}^\mu$ in terms of the canonical basis.

\begin{corollary}
\label{cor:multiplicity through BB}
The coefficient $c_{\lambda, \nu}^\mu$ is equal to the number of vectors
$b \in \BB$ of degree $\lambda + \nu - \mu$ satisfying the following property:
if $b = b_\ii (t_1, \dots, t_m)$, and $\ii \in R(w_0)$ begins with $i$ and ends with $j$ then
$t_1 \leq  \langle \alpha_i^\vee, \lambda\rangle$ and
$t_m \leq \langle \alpha_{j^*}^\vee, \nu \rangle$.
\end{corollary}

\subsection{Transition maps}
\label{sec:transition}
Let $\ii$ and $\ii'$ be two reduced words for $w_0$.
In view of Proposition~\ref{pr:characterization of BB},
there is a bijective \emph{transition map}
$$R_\ii^{\ii'} =
(b_{\ii'})^{-1} \circ b_\ii:\ZZ_{\ge 0}^m \to \ZZ_{\ge 0}^m $$
between the corresponding Lusztig parametrizations of the canonical basis $\BB$.
It turns out that each component of a tuple $R_\ii^{\ii'} (t)$
can be expressed through the components of $t$
as a ``tropical" subtraction-free rational expression.

\begin{example}
\label{ex:A2 transition}
{\rm Let $\gg = sl_3$ be of type $A_2$, and let $\ii = (1,2,1)$ and $\ii' = (2,1,2)$
be the two reduced words for $w_0$.
The transition map $R_\ii^{\ii'}$ between two Lusztig parametrizations
was computed in \cite{lu}: the components of $t' = R_\ii^{\ii'}(t)$ are given by
$$t'_1 = t_2 + t_3 - \min \ (t_1, t_3), \,\,  t'_2 = \min \ (t_1, t_3), \,\,
t'_3 = t_1 + t_2 - \min \ (t_1, t_3) \ ,$$
which can also be written as
\begin{equation}
\label{eq:Lustzig A2 transition}
t'_1 = \left[\frac{t_2 t_3}{t_1 + t_3}\right]_{\rm trop}, \,\,
t'_2 = [t_1 + t_3]_{\rm trop}, \,\,
t'_3 = \left[\frac{t_1 t_2}{t_1 + t_3}\right]_{\rm trop} \ .
\end{equation}
}
\end{example}

Lusztig observed that the same formula (\ref{eq:Lustzig A2 transition})
understood in a usual (non-tropical) sense describes the transition map
that relates the parameters in two factorizations of a matrix
$x \in N \subset G = SL_3$: if
$$x = x_1 (t_1) x_2 (t_2) x_1 (t_3) = x_2 (t'_1) x_1 (t'_2) x_2(t'_3)$$
then the $t'_i$ are expressed through the $t_i$ via
$$t'_1 = \frac{t_2 t_3}{t_1 + t_3}, \,\,
t'_2 = t_1 + t_3, \,\,
t'_3 = \frac{t_1 t_2}{t_1 + t_3} \ . $$
The birational transformation $(t_1,t_2,t_3) \mapsto (t'_1, t'_2, t'_3)$
becomes a bijection if we restrict it to
the tuples $(t_i)$ and $(t'_i)$ of \emph{positive} real numbers.
To generalize this observation to an arbitrary semisimple group $G$,
we will need the notion of \emph{total positivity}.

\subsection{Total positivity in $N$}
Following G.~Lusztig \cite{lusztig-reductive}, we define the subset $G_{\geq 0}$ of
\emph{totally nonnegative} elements in $G$
as the multiplicative monoid with unit generated by the elements
$t^{\alpha_i^\vee}$, $x_i (t)$, and $y_i (t)$ for all $i$ and all positive real $t$.
As proved in \cite{fzproc}, $x \in G$ is totally nonnegative if and only if
all generalized minors $\Delta_{\gamma, \delta}$ take nonnegative real values at $x$.
In these notes, we shall only need the set $N_{>0}$ of
\emph{totally positive} elements of $N$.
It can be defined as follows:
\begin{equation}
\label{eq:tp-N}
N_{> 0} = G_{\geq 0} \cap N \cap B_- w_0 B_- \ .
\end{equation}
The following proposition which is an easy consequence
of results in \cite{lusztig-reductive,fz} provides three more equivalent
definitions of $N_{> 0}$.

\begin{proposition}
\label{pr:tp-N}
An element $x \in N$ is totally positive if and only if it satisfies
each of the following conditions:

\begin{itemize}

\item [{\rm (1)}] $\Delta_{\gamma, \delta} (x) > 0$ for any upper-triangular
$(\gamma, \delta)$.

\item [{\rm (2)}] For any reduced word $\ii = (i_1, \dots, i_m) \in R(w_0)$,
there is a unique factorization $x = x_{i_1} (t_1) \cdots  x_{i_m} (t_m)$
with all $t_k$ positive real numbers.

\item [{\rm (3)}] $x \in G_{\geq 0}$, and
$\Delta_{\omega_i, w_0 \omega_i} (x) > 0$ for all $i$.

\end{itemize}
\end{proposition}

In view of condition (2) in Proposition~\ref{pr:tp-N},
any reduced word $\ii = (i_1, \dots, i_m) \in R(w_0)$
gives rise to a bijection $x_\ii: \RR_{>0}^m \to N_{> 0}$
given by
$$x_\ii (t_1, \dots, t_m) = x_{i_1} (t_1) \cdots  x_{i_m} (t_m) \ .$$
It follows that, for any $\ii, \ii' \in R(w_0)$,
there is a bijective transition map
$$\tilde R_\ii^{\ii'} =
(x_{\ii'})^{-1} \circ x_\ii:\RR_{> 0}^m \to \RR_{> 0}^m $$
that relates the corresponding parametrizations of
the totally positive variety $N_{>0}$.
We shall use the notation $\bigl(\tilde R_\ii^{\ii'}\bigr)^\vee$ for the transition maps
defined in the same way for the group $G^\vee$ associated with the Langlands dual
Lie algebra $\gg^\vee$.
The following theorem is a special case of \cite[Theorem~5.2]{bz01}.

\begin{theorem}
\label{th:geometric lifting}

\noindent

\begin{itemize}

\item[{\rm (i)}] 
Each component of
$\tilde R_\ii^{\ii'} (t)$ is a subtraction-free rational expression in the components
of $t$.

\item[{\rm (ii)}] 
Each component of
Lusztig's transition map $R_\ii^{\ii'}$ between two parametrizations
of the canonical basis is the tropicalization of the corresponding component of
$\bigl(\tilde R_\ii^{\ii'}\bigr)^\vee$.
\end{itemize}
\end{theorem}

\emph{Outline of the proof.}
We use the following well known result due to Tits:
every two reduced words of the same element of $W$ can be obtained from each other
by a sequence of braid moves.
Using this result, we reduce both claims in Theorem~\ref{th:geometric lifting}
to the rank 2 case.
The latter is done by a direct case-by-case computation.

\subsection{Completing the proof of Theorem~\ref{th:LR-Plucker-Lusztig}}
Let us consider the variety $\MM$ introduced in Section~\ref{sec:main}.
Let $\MM (\RR_{\geq 0})$ denote the ``totally positive" part of $\MM$ formed
by the tuples $(M_{\gamma, \delta})$ whose components are positive real numbers.
The following proposition is a special case of \cite[Theorem~5.13]{bz01}.

\begin{proposition}
\label{pr:N-M}
The map $B \to \MM$ given by
$x \mapsto (M_{\gamma, \delta} = \Delta_{\gamma, \delta}(x))$
restricts to a bijection between $N_{> 0}$ and
the set $\MM_{> 0}^{\rm unip}$ of all tuples
$(M_{\omega_i,\gamma}) \in  \MM_{> 0}$ such that
$M_{\omega_i,\omega_i} = 1$ for all $i$.
\end{proposition}

Now let us fix a reduced word $\ii \in R(w_0)$.
Combining the bijection $N_{> 0} \to \MM_{> 0}^{\rm unip}$
in Proposition~\ref{pr:N-M}
with the bijection $x_\ii: \RR_{>0}^m \to N_{> 0}$,
we obtain a bijection $\mu_\ii: \RR_{>0}^m \to \MM_{> 0}^{\rm unip}$.
Both $\mu_\ii$ and its inverse are given by subtraction-free rational expressions
but nice explicit expressions for them are not known.
However, we shall only need some partial information: if
$\mu_\ii (t_1, \dots, t_m) = (M_{\gamma, \delta})$ then we have
\begin{equation}
\label{eq:t1-tm}
t_1 = \frac{M_{\omega_{i_1}, w_0 \omega_{i_1}}}
{M_{s_{i_1} \omega_{i_1}, w_0 \omega_{i_1}}} \ , \,\,
t_m = \frac{M_{\omega_{i_m^*}, w_0 \omega_{i_m^*}}}
{M_{\omega_{i_m^*}, s_{i_m} w_0 \omega_{i_m^*}}}
\end{equation}
(cf. formula (4.19) in \cite{bz01}), and
\begin{equation}
\label{eq:special-minors}
M_{\omega_i, s_i \omega_i} = \sum_{k: i_k = i} t_k \ .
\end{equation}

Now all the ingredients are in place for completing the proofs of
Proposition~\ref{pr:Un-basis} and Theorem~\ref{th:LR-Plucker-Lusztig}.
With some abuse of notation, let us denote by the same symbol $\BB$
the canonical basis of $U^+$ and the basis of $U(\nn)$ obtained
from it by specializing $q$ to $1$.
For the latter basis, property (1) in Proposition~\ref{pr:Un-basis} follows from
Proposition~\ref{pr:Lusztig parametrizations}.
Combining the inverse of Lusztig's parametrization $b_\ii$ with
the tropical version of the bijection $\mu_\ii^\vee$ (here
$\mu_\ii^\vee$ is defined in the same way as $\mu_\ii$ above but
for the Langlands dual Lie algebra), we obtain an embedding
$\BB \to \MM_{\rm trop}^\vee$.
By part (ii) of Theorem~\ref{th:geometric lifting}, this embedding
does not depend on the choice of a reduced word $\ii$.
It remains to show that the image of this embedding is given by the
conditions in Theorem~\ref{th:LR-Plucker-Lusztig}.
Now condition $(0)$ is the tropical version of the equality
$M_{\omega_i,\omega_i} = 1$ in
Proposition~\ref{pr:N-M}.
By taking the tropical version of (\ref{eq:special-minors}),
we see that condition $(1)$ simply means that all components
$t_k$ in Lusztig's parametrization of the canonical basis are
nonnegative.
Finally, conditions $(2)-(4)$ are the result of rewriting the conditions
in Corollary~\ref{cor:multiplicity through BB} in terms of the
$M_{\gamma, \delta}$; this is done with the help of
(\ref{eq:t1-tm}).

\begin{remark}
{\rm The inverse bijection $x_\ii^{-1}: N_{> 0} \to \RR_{> 0}^m$ was computed
in \cite{bz97} (for the type $A_r$, this was done in \cite{bfz}).
Further generalizations of this result can be found in \cite{fz,bz01}.
An interesting feature of the answer is that the parameters
$t_k$ in the factorization $x = x_\ii (t_1, \dots, t_m)$ are expressed
(in a quite simple way) in terms of generalized minors evaluated not at $x$
but at another element obtained from $x$ by some birational transformation
(``twist").
However, in view of Corollary~\ref{cor:multiplicity through BB}, to compute
the multiplicities $c_{\lambda \nu}^\mu$, we only need explicit expressions
for the first and the last of the $t_k$; luckily, these two parameters can be
expressed in terms of minors of $x$ via (\ref{eq:t1-tm}).
}
\end{remark}

\section{Lecture III: Introduction to cluster algebras}
\label{sec:lecture3}

\subsection{Motivations, main features and examples}
After providing the canonical basis $\BB$ with a combinatorial
parametrization, the next natural step is to study the algebraic structure of $\BB$.
We shall focus on the classical case $q=1$, and, instead of talking about
the basis $\BB$ in $U(\nn)$, turn our attention to the dual basis
$\BB^{\rm dual}$ in the ring of regular functions $\CC[N]$
(recall from Section~\ref{sec:dual} that each homogeneous component
$\CC[N](\gamma)$ is naturally identified with the dual space~$(U(\nn)(\gamma))^*$).
One can expect $\BB^{\rm dual}$ to be more tractable since, in contrast to $U(\nn)$,
the algebra $\CC[N]$ is commutative.
(This distinction disappears under the $q$-deformation: the algebra
$\CC_q [N]$ is not only non-commutative but in fact isomorphic to
$U^+ = U_q(\nn)$.)
Another motivation for studying $\BB^{\rm dual}$ comes from the
theory of total positivity: this connection is given by an important result
(essentially due to G.~Lusztig \cite{lusztig-reductive}) that every function from
$\BB^{\rm dual}$ takes positive real values on
the totally positive variety $N_{> 0}$.

The dual canonical basis $\BB^{\rm dual}$ was constructed
explicitly in several small rank cases: for the types
$A_2$, $B_2$ and $A_3$ this was done, respectively in \cite{gz86},
\cite{rz} and \cite{bz1}.
In all these cases, $\BB^{\rm dual}$ consists of certain
monomials in a distinguished family of generators.
(A word of warning: this basis is very different from the standard
monomial basis studied extensively by Lakshmibai, Seshadri and
their collaborators.)
The monomials that constitute $\BB^{\rm dual}$ are defined by not allowing certain pairs of
generators to appear together.
In each case, the product of every two ``incompatible"  generators
can be expressed as the sum of two allowed monomials.
Such \emph{binomial exchange relations} turn out to play a very
important part in the study of canonical
bases and total positivity undertaken by the author and his collaborators
(see, e.g., \cite{bfz,bz97,bz01,fz,ssvz,z-imrn}); important
instances of these relations are presented in
Propositions~\ref{pr:minors-Dodgson} and \ref{pr:minors-Plucker}.

Inspired by all this work,
a new class of commutative algebras called \emph{cluster algebras}
was introduced in \cite{fzclus1} as an attempt to create
a natural algebraic framework for the study of the dual canonical
basis and total positivity.
Before discussing precise definitions, let us present some of the main features
of cluster algebras.

For any positive integer $n$,
a cluster algebra $\AA$ of rank $n$ is a commutative algebra
with unit and no zero divisors over a finite polynomial ring,
equipped with a distinguished family of generators
called \emph{cluster variables}. The set of cluster variables is
the (non-disjoint) union of a distinguished collection of
$n$-subsets called \emph{clusters}.
These clusters have the
following \emph{exchange property}: for any cluster $\xx$ and
any element $x \in \xx$, there is another cluster obtained from
$\xx$ by replacing $x$ with an element $x'$ related to $x$ by a
binomial exchange relation
$$x x' = M_1 + M_2 \ ,$$
where $M_1$ and $M_2$ are disjoint monomials in
the $n-1$ variables $\xx - \{x\}$. Furthermore, any two
clusters can be obtained from each other by a sequence of
exchanges of this kind.

The prototypical example of a cluster algebra of rank 1
is the coordinate ring $\AA = \CC [SL_2]$ of the group $SL_2$, viewed
in the following way.
Writing a generic element of $SL_2$ as $\mat{a}{b}{c}{d},$
we consider the entries $a$ and $d$ as cluster variables,
and the entries $b$ and $c$ as scalars.
There are two clusters $\{a\}$ and $\{d\}$,
and $\AA$ is the algebra over the ground ring
$\CC [b,c]$ generated by the cluster variables $a$ and $d$
subject to the binomial exchange relation
$$ad = 1 + b c \ .$$

Another important incarnation of a cluster algebra of rank 1
is the coordinate ring $\AA = \CC [N_- \backslash G]$ of the base affine space
of the special linear group $G = SL_3$; here $N$ is the maximal unipotent subgroup
of $SL_3$ consisting of all unipotent upper triangular matrices.
Using the standard notation $(x_1, x_2, x_3, x_{12}, x_{13},
x_{23})$ for the Pl\"ucker coordinates on $N_- \backslash G$, we view $x_2$ and $x_{13}$
as cluster variables;
then $\AA$ is the algebra over the polynomial ring
$\CC [x_1,x_3, x_{12}, x_{23}]$ generated by the two cluster variables
$x_2$ and $x_{13}$
subject to the binomial exchange relation
$$x_2 x_{13} = x_1 x_{23} + x_3 x_{12} \ .$$
(Note that in this example the dual canonical basis consists of all monomials in
the six Pl\"ucker coordinates that are not divisible by $x_2 x_{13}$.)

A nice example of a cluster
algebra of an arbitrary rank $n$ is the homogeneous coordinate ring
$\CC [Gr_{2,n+3}]$ of the Grassmannian of $2$-dimensional
subspaces in~$\CC^{n+3}$.
This ring is generated by the Pl\"ucker coordinates $x_{ij}$,
for $1 \leq i < j \leq n+3$, subject to the relations
$$x_{ik} x_{jl} = x_{ij} x_{kl} + x_{il} x_{jk} \ ,$$
for all $i < j < k < l$.
It is convenient to identify the indices $1, \dots, n+3$ with the
vertices of a convex $(n+3)$-gon, and
the Pl\"ucker coordinates with its sides and diagonals.
Let us view the sides of the polygon as
scalars, and the diagonals as cluster variables.
The clusters are the maximal families of pairwise non-crossing
diagonals; thus, they are in a natural bijection with the
triangulations of the polygon.
It is known that the monomials in Pl\"ucker coordinates not involving crossing diagonals
form a linear basis in~$\CC [Gr_{2,n+3}]$.
To be more specific, this ring is naturally identified with the ring of
polynomial $SL_2$-invariants of an $(n+3)$-tuple of points in~$\CC^2$.
Under this isomorphism, the above basis corresponds to the basis
considered in \cite{kungrota84,stur}.

\subsection{Cluster algebras of geometric type}
For the sake of simplicity, we shall only
discuss a special
class of cluster algebras from \cite{fzclus1}
called \emph{cluster algebras of geometric type};
this case covers most applications that we have in mind.

Let $I \subset \tilde I$ be two finite index sets of cardinalities
$|I| = n \leq m = |\tilde I|$.
We also denote $I^c = \tilde I - I$, and introduce a family of
commuting variables $p_i$ for $i \in I^c$.
Let $\tilde B = (b_{ij})$ be an integer $m \times n$ matrix with
rows indexed by $\tilde I$ and columns indexed by $I$.
We call the submatrix $B = (b_{ij})_{i,j \in I}$ the
\emph{principal part} of $\tilde B$, and the complementary
submatrix $B^c = (b_{ij})_{i \in I^c, j \in I}$ the
\emph{complementary part}.
The only requirement on $\tilde B$ is that its principal part $B$
must be \emph{skew-symmetrizable}, i.e., $DB$ is skew-symmetric for
some diagonal $I \times I$ matrix $D$ with positive integer
diagonal entries.
To every such $\tilde B$, we associate the cluster algebra
$\AA (\tilde B)$ over the ground ring $S = \ZZ[p_i: i \in I^c]$.
This is done in four steps.

\smallskip

\noindent \emph{Step 1.}
Let $\TT_n$ denote the $n$-regular tree whose edges are
labeled by the elements of $I$
so that the $n$ edges emanating from each vertex receive
different labels.
We will write $t \overunder{i}{} t'$ if vertices
$t,t'\in\TT_n$ are joined by an edge labeled by~$i$.
To every vertex $t \in \TT_n$, we associate a {\it cluster}
$\xx (t)$ consisting of $n$ independent \emph{cluster variables}
$x_i (t)$ for $i \in I$.
We also set $x_i (t) = p_i$ for $i \in I^c$ and $t \in \TT_n$.

\smallskip

\noindent \emph{Step 2.}
Pick a vertex $t_0 \in \TT_n$, and set $\tilde B (t_0) = \tilde B$.
Using the matrix $\tilde B$ as a starting point, let us now construct a
family of integer $\tilde I \times I$-matrices $(\tilde B (t))$, one for each vertex
$t \in \TT_n$.
This is done with the help of the following operation called
\emph{matrix mutation}: for
every $k \in I$, define ${\tilde B}' = \mu_k (\tilde B)$
by setting
\begin{equation}
\label{eq:mutation}
b'_{ij} =
\begin{cases}
-b_{ij} & \text{if $i=k$ or $j=k$;} \\[.15in]
b_{ij} + \displaystyle\frac{|b_{ik}| b_{kj} +
b_{ik} |b_{kj}|}{2} & \text{otherwise.}
\end{cases}
\end{equation}
The following properties are checked immediately.

\begin{proposition}
\label{pr:matrix-mutation}

\noindent
\begin{itemize}
\item[{\rm (i)}] The operation $\mu_k$ is an involution.

\item[{\rm (ii)}] The principal part of
$\mu_k (\tilde B)$ is $\mu_k (B)$, where
$B$ is the principal part of~$\tilde B$.

\item [{\rm (iii)}] If $B$ is skew-symmetrizable then so is $\mu_k (B)$.
\end{itemize}
\end{proposition}

The family of matrices $(\tilde B (t))$ is now uniquely determined
by the following requirement: $\tilde B(t') = \mu_k (\tilde B(t))$
whenever $t \overunder{k}{} t'$.

\smallskip

\noindent \emph{Step 3.}
The binomial exchange relations for cluster variables take the
following form: for every edge $t \overunder{j}{} t'$ in $\TT_n$,
we write
\begin{equation}
\label{eq:exchange1}
x_i (t') = x_i (t) \,\, (i \neq j)
\end{equation}
and

\begin{equation}
\label{eq:exchange2}
x_j (t) x_j (t') = \prod_{i \in \tilde I} x_i (t)^{\max (b_{ij} (t),0)} +
\prod_{i \in \tilde I} x_i (t)^{\max (- b_{ij} (t),0)}.
\end{equation}
Note that, in view of (\ref{eq:mutation}) and (\ref{eq:exchange1}),
the relation (\ref{eq:exchange2}) is symmetric in $t$ and $t'$.

\smallskip

\noindent \emph{Step 4.}
Moving away from $t_0$ and iterating (\ref{eq:exchange1}) and (\ref{eq:exchange2})
along the way, we can express every cluster variable $x_i (t)$
as a rational function in the initial
cluster $\xx (t_0)$ with coefficients in the ground ring $S$.
The cluster algebra $\AA(\tilde B)$ is defined as the
$S$-subalgebra of the field of rational functions $S(\xx(t_0))$
generated by all cluster variables.
Note that this definition does not depend on the choice of an
initial cluster.
Furthermore, $\AA(\tilde B)$ depends only on the
mutation-equivalence class of a matrix $\tilde B$,
i.e., $\AA(\tilde B)$ is naturally identified with $\AA(\tilde B')$
for any $\tilde B'$ obtained from $\tilde B$ by a sequence of
matrix mutations.

\subsection{The Laurent phenomenon}
One of the main structural features of cluster algebras established in
\cite{fzclus1} is the following \emph{Laurent phenomenon}.

\begin{theorem}
\label{th:cluster-Laurent}
Any cluster variable viewed as a  rational function in the variables of any given cluster
is a Laurent polynomial whose coefficients are integer Laurent polynomials in the variables $p_i$.
\end{theorem}

This is quite surprising: the numerators of these Laurent polynomials
may contain a huge number of monomials, and the numerator for a cluster variable $x$
moves into the denominator when we compute the cluster variable $x'$ obtained
from $x$ by an exchange (\ref{eq:exchange2}).
The magic of the Laurent phenomenon is that, at every stage of
this recursive process, a cancellation will inevitably occur,
leaving a single monomial in the denominator.

As a byproduct of this development, a number of other instances
of the Laurent phenomenon spreading beyond the cluster algebra framework
were established in \cite{fzlaurent}.
The results include a proof of the conjecture
of D.~Gale and R.~Robinson on integrality of generalized Somos sequences, as well as proofs of the
Laurent property for several multidimensional recurrences,
conjectured by J.~Propp, N.~Elkies, and M.~Kleber.

\subsection{Classification of cluster algebras of finite type}
We say that a cluster algebra is of \emph{finite type} if it has finitely many
distinct cluster variables (or equivalently, finitely many
distinct clusters).
A classification of such algebras will be given in a forthcoming
paper by S.~Fomin and the author.
Remarkably, this classification provides yet another instance of
the celebrated Cartan-Killing classification.
To be more specific, to every skew-symmetrizable matrix $B=(b_{ij})$
one can associate a matrix $A = A(B)=(a_{ij})$ of the same size
by setting
\begin{equation}
\label{eq:cartan matrix}
a_{ij} =
\begin{cases}
2 & \text{if $i=j$;} \\[.1in]
- |b_{ij}| & \text{if $i\neq j$.}
\end{cases}
\end{equation}
The matrix $A$ satisfies the first two conditions in
Definition~\ref{def:Cartan-matrix}; such matrices are often
referred to as generalized Cartan matrices.

\begin{theorem}
\label{th:CK-classification}
Let $\tilde B$ be an integer $m \times n$ matrix
whose principal part $B$ is skew-symmetrizable.
The cluster algebra $\AA (\tilde B)$ is of finite type if and only
if $B$ is mutation equivalent to a matrix $B'$ such that
$A(B')$ is a Cartan matrix of finite type, i.e., is a direct sum of the matrices
from the Cartan-Killing list $A_n, B_n, \dots, G_2$.
\end{theorem}

For cluster algebras of rank $2$, this was proved in \cite{fzclus1}.
The algebra $\CC[Gr_{2,n+3}]$ discussed above provides an example
of a cluster algebra of finite type; its Cartan-Killing type is $A_n$.

We expect that, for any simply-connected connected semisimple group
$G$, each of the coordinate rings $\CC [G]$, $\CC [N_-\backslash G]$,
$\CC [N]$, as well as
coordinate rings of many other interesting varieties related
to~$G$, have a natural structure of a cluster algebra.
By Theorem~\ref{th:CK-classification}, those among these
cluster algebras that are of finite type have some kind of
hidden symmetry associated to a Cartan-Killing type that may be
completely different from the Cartan-Killing type of the ambient
group $G$.
Here are some instances of this phenomenon (in each case, a
cluster algebra is accompanied by its ``cluster type'' given by
Theorem~\ref{th:CK-classification}):
\begin{center}
\begin{tabular}{ccccc}
$\CC[{\rm Gr}_{2,n+3}]$ & $A_n$ &\qquad& $\CC[N_-\backslash SL_3]$ & $A_1$\\[.1in]
$\CC[{\rm Gr}_{3,6}]$   & $D_4$ && $\CC[N_-\backslash SL_4]$ & $A_3$ \\[.1in]
$\CC[{\rm Gr}_{3,7}]$   & $E_6$ && $\CC[N_-\backslash SL_5]$ & $D_6$ \\[.1in]
$\CC[{\rm Gr}_{3,8}]$   & $E_8$ && $\CC[N_-\backslash Sp_4]$ & $B_2.$
\end{tabular}
\end{center}

Theorem~\ref{th:CK-classification} suggests that there must be
some connection between a cluster algebra $\AA$ of finite type
and the semisimple Lie algebra $\gg$  associated with its cluster
type. Here is just one such connection (it will be proved among many
others in a forthcoming sequel to \cite{fzclus1}): the number of
distinct cluster variables in $\AA$ turns out to be equal to the
dimension of the base affine space for $\gg$ (equivalently, this
is the rank of $\gg$ plus the number of positive roots).
For instance, since the cluster algebra $\CC[N_-\backslash SL_5]$
has cluster type $D_6$, it has $36$ cluster variables.

We naturally expect that even when a cluster algebra $\AA(\tilde B)$
is of infinite type, it is closely connected with the
Kac-Moody algebra associated to the generalized Cartan matrix $A(B)$.
Although some promising first steps in this direction were made in
\cite{fzclus1}) for the rank $2$ case,
the general case is still far from being understood.

\textsc{Acknowledgements.}
I thank the organizers of the NATO Advanced Study Institute
for the invitation to give these lectures, and
the Isaac Newton Institute for Mathematical
Sciences for support and hospitality.
Special thanks to
Sergey Fomin whose hard work helped to make this meeting truly enjoyable.
I am also grateful to Sergey for helpful editorial suggestions.

\end{document}